\documentclass[12pt]{article}

\usepackage{fullpage, verbatim, xy, tikz-cd}
\usepackage{amsmath,amssymb,latexsym, amsfonts, amscd, amsthm, mathrsfs}
\xyoption{all}
\usepackage[none]{hyphenat}
 \usepackage[bibtex-style]{amsrefs}

\renewcommand{\tilde}{\widetilde}

\DeclareMathOperator{\Stab}{Stab}

\DeclareMathOperator{\BK}{BK}
\DeclareMathOperator{\BKd}{BKd}
\DeclareMathOperator{\Vol}{Vol}

\DeclareMathOperator{\reg}{reg}
\DeclareMathOperator{\gr}{gr}
\DeclareMathOperator{\pr}{pr}

\DeclareMathOperator{\Span}{Span}

\newcommand{\sL}{\mathscr{L}}

\DeclareMathOperator{\End}{End}

\DeclareMathOperator{\Gal}{Gal}

\DeclareMathOperator{\GL}{GL}

\DeclareMathOperator{\Ind}{Ind}

\DeclareMathOperator{\rank}{rank}

\def\N{\mathrm{N}}

\newcommand{\Z}{\mathbf{Z}}

\newcommand{\Q}{\mathbf{Q}}
\newcommand{\C}{\mathbf{C}}

\newcommand{\R}{\mathbf{R}}

\newcommand{\cO}{\mathcal{O}}

\newcommand{\fm}{\mathfrak{m}}

\newcommand{\fp}{\mathfrak{p}}

\newtheorem{lemma}{Lemma}[section]

\newtheorem{prop}[lemma]{Proposition}

\newtheorem{theorem}[lemma]{Theorem}
\newtheorem{conj}[lemma]{Conjecture}
\newtheorem{prob}[lemma]{Problem}

\theoremstyle{remark}

\theoremstyle{definition}
\newtheorem{definition}[lemma]{Definition}

\newtheorem{remark}[lemma]{Remark}
\newtheorem{question}[lemma]{Question}

\title{Ranks of Matrices of Logarithms of Algebraic Numbers II: The Matrix Coefficient Conjecture}

\author{Samit Dasgupta \\ Mahesh Kakde}

\begin{document}
\baselineskip 15.8pt

\makeatletter
\renewcommand{\maketitle}{\bgroup\setlength{\parindent}{0pt}
\begin{center}
  \Large{\textbf{\@title}}
  
  \vspace{4mm}
  
  \large{\textsc{\@author}}
  
  \vspace{4mm}
\end{center}\egroup
}
\makeatother

\maketitle

{\small 
\paragraph{Abstract.}
Many questions in number theory concern the nonvanishing of determinants of square matrices of logarithms (complex or $p$-adic) of algebraic numbers.  We present a new conjecture that states that if such a matrix has vanishing determinant, then after a rational change of basis on the left and right, it can be made to have a vanishing coefficient.
}

 \tableofcontents

\section{Introduction}

There are many settings in number theory where one wants to know that certain square matrices of logarithms (classical or $p$-adic) have non-vanishing determinant.
Leopoldt's conjecture and the Gross--Kuz'min conjecture are two well-known examples (these conjectures are stated precisely in the proofs of Theorems~\ref{t:leopoldt} and~\ref{t:gk} below).  Both conjectures are special cases of 
the Structural Rank Conjecture of transcendence theory, whose statement we now recall.  To streamline the discussion we discuss the archimedean and $p$-adic contexts simultaneously.
Let \[ \sL = \{x \in \C \colon e^x \in \overline{\Q} \}  \] 
denote the $\Q$-vector space of logarithms of algebraic numbers.
Similarly, for each prime $p$ we fix an embedding $\overline{\Q} \subset \C_p$ and let
\[ \sL_p = \{ \log_p(x) \colon x \in \overline{\Q}^*\} \subset \C_p, \]
where $\log_p$ denotes Iwasawa's $p$-adic logarithm.
\begin{conj}[Structural Rank Conjecture] \label{c:rank} Let $M$ be an $m \times n$ matrix with coefficients in $\sL$ (or $\sL_p$). Write \[ M = \sum_{i=1}^{k} M_i \lambda_i \] where $M_i \in M_{m \times n}(\Q)$ and the $\lambda_i$ are elements of $\sL$ (or $\sL_p$) that are linearly independent over ${\Q}$.  Then the rank of $M$ is equal to the rank of 
 $ \sum_{i=1}^{k} M_i x_i$, viewed as a matrix with coefficients in the field $\Q(x_1, \dots, x_k)$.
\end{conj}

We discussed the Structural Rank Conjecture in our expository work \cite{das};  for example, we recalled Roy's proof \cite{roy} of the equivalence of the Structural Rank Conjecture to a special case of Schanuel's conjecture.
The goal of the present article is to state and explore a new conjecture that is ostensibly weaker than the Structural Rank Conjecture, but is still strong enough to imply important open problems in transcendence theory, including the Four Exponentials conjecture, Leopoldt's conjecture, and the Gross-Kuz'min conjecture.
We call this conjecture the Matrix Coefficient Conjecture because it states that after a rational change of bases on the left and right, a singular square matrix of logarithms of algebraic numbers can be made to have some coefficient vanish.  In the statement below, $\langle \ , \ \!\rangle$ denotes the usual dot product.

\begin{conj}[Matrix Coefficient Conjecture]  \label{c:mainq} Let $M \in M_{n \times n}(\sL)$ or $M_{n \times n}(\sL_p)$ and suppose that $\det(M) = 0$.  Then there exist nonzero vectors $w, v \in \Q^n$ such that $\langle w, Mv \rangle = 0$.
\end{conj}

It is easy to show that in the case $n=2$, the Matrix Coefficient Conjecture is equivalent to the Four Exponentials Conjecture.
The main contributions of the paper are as follows:

\begin{itemize}
\item We prove that the Structural Rank Conjecture implies the Matrix Coefficient Conjecture.
\item We prove that the Matrix Coefficient Conjecture ($\sL_p$ version) implies Leopoldt's conjecture and the Gross--Kuz'min conjecture.
\item We propose a strategy for attacking the  Matrix Coefficient Conjecture modeled on the theorem of Waldschmidt--Masser that gives a partial result toward the 
Structural Rank Conjecture.  Importantly, we replace the condition on degree of auxiliary polynomials in Masser's result with a condition on the support of the polynomials (i.e. the set of exponents of monomials appearing with nonzero coefficients).  We prove some of the implications necessary in our strategy, leaving open the implication written ``(w) $\Longrightarrow$ (m)'' below that would be needed to complete the proof.
\end{itemize}

We would like to thank Damien Roy and Michel Waldschmidt for helpful discussions regarding an earlier draft of this article.  We are indebted to Ezra Miller, J.~Maurice Rojas, Henry Schenck, and Frank Sottile for helpful discussions regarding \S\ref{s:strategy},
 in particular to Prof.~Rojas for suggesting the idea of appending $n-i$ copies of a generic linear polynomial $\ell$. The first named author is supported by NSF grant DMS--2200787. The second named author is supported by DST-SERB grant SB/SJF/2020-21/11, Tata Education and Development Trust, and DST FIST program - 2021 [TPN - 700661].

\section{Relation with the Structural Rank Conjecture}

The following is a comforting consistency check.

\begin{theorem}  \label{t:imply} Conjecture~\ref{c:rank} implies Conjecture~\ref{c:mainq}
\end{theorem}

This follows from the following.

\begin{theorem} \label{t:singular}  Let $F$ be an infinite field and let $V$ be a finite dimensional vector space over $F$.  Let $E \subset \End(V)$ be a subspace such that every element of $E$ is singular.  Then there exist nonzero $v \in V, w \in V^*$ such that $\langle w, Av \rangle = 0$ for all $A \in E$.
\end{theorem}

Theorem~\ref{t:imply} immediately follows from Theorem~\ref{t:singular} by letting $E$ be the $\Q$-span of the $M_i$. 
We thank Damien Roy for communicating the following proof to us; our original proof was more complicated than what is written below.

\begin{proof}
We induct on $n = \dim V$.  If $V$ is 1-dimensional over $F$, the result is trivial.
In the general case, fix a nonzero element $A \in E$.  After applying an element of $\GL(V)$ on the left and right---which changes neither the assumptions nor the conclusion of the theorem---we can assume that $A$ is a diagonal matrix with first $r$ diagonal entries equal to $1$ and last $n-r$ diagonal entries equal to 0, where $1 < r < n$.  Note that $r \neq 0$ since $A \neq 0$ and $r \neq n$ since $\det(A)=0$.

We now evaluate $\det(At + B)$ where we view $t$ as an indeterminate and $B \in E$ also as a generic variable.
The leading coefficient of this expression as a polynomial in $t$ is $\det(B')$, where $B'$ is the lower right $(n-r) \times (n-r)$ block of $B$.  Since $F$ is infinite, and this polynomial is by assumption 0 for all values of $t \in F$, it follows that $\det(B') = 0$ for all $B \in E$.  By induction, there exist nonzero $w', v' \in F^{n-r}$ such that $\langle w', B' v \rangle = 0$ for all $B \in E$.  Padding $w'$ and $v'$ with zeroes we obtain $w, v \in V = F^n$ such that that $\langle w, B v \rangle = 0$ for all $B \in E$ as desired.
\end{proof}

It is natural to formulate the following $\overline{\Q}$-version of Conjecture~\ref{c:mainq}.

\begin{conj} \label{c:main}  Let $M$ be an $n \times n$ matrix  whose elements lie in the $\overline{\Q}$ vector space spanned by $\sL$ (or $\sL_p$). Suppose that $\det(M) = 0$.  Then there exist nonzero vectors $w, v \in \overline{\Q}^n$ such that $\langle w, Mv \rangle = 0$.
\end{conj}

 Roy proved that a suitably strong version of the Waldschmidt--Masser result over $\Q$ implies the stronger version over $\overline{\Q}$.  This leads naturally to the following problem.

\begin{prob} Show that Conjecture~\ref{c:mainq} or a suitable strengthening that remains ``over $\Q$" implies Conjecture~\ref{c:main}.
\end{prob}

 Roy proved that the Waldschmidt--Masser result can be generalized to $\Q+ \sL$, the $\Q$-vector space spanned by $\Q$ and $\sL$, and to the $\overline{\Q}$ vector space spanned by $\Q + \sL$ (or $\sL_p$ in both instances).  Since Conjecture~\ref{c:rank} is expected to hold with $\sL$ replaced by $\Q + \sL$, 
the direct generalization of Conjecture~\ref{c:mainq} to  $\Q+ \sL$ should hold as well; we thank Roy for pointing this out to us.

\section{Connection to Regulators in Number Theory}

\begin{theorem} Conjecture~\ref{c:mainq} implies Leopoldt's conjecture. \label{t:leopoldt}
\end{theorem}

\begin{proof}  It is well-known that it suffices to consider  finite Galois extensions $F/\Q.$  Write \[ G = \Gal(F/\Q) \] and fix once and for all embeddings $F \subset \overline{\Q} \subset \C_p$.  
Let $U = \cO_F^* \otimes_{\Z} \Q$.
Write 
\[ n = \#G = [F:\Q], \qquad r = \dim_{\Q} U = 
 \begin{cases} n-1 & F \text{ is real} \\
n/2 - 1 & F \text{ is complex.} \end{cases}  \]
Leopoldt's conjecture is that the map 
\[ \reg_p \colon U \otimes_\Q \C_p \longrightarrow \bigoplus_{\tau \in G} \C_p, \qquad u \mapsto (\log_p \tau^{-1}(u)) \]
has image of dimension $r$.  The map $\reg_p$ is $G$-equivariant, and we can describe the $G$-module structure of $U$ explicitly.  Namely, if $\infty$ denotes some fixed archimedean place of $F$ and $G_\infty \subset G$ is the associated decomposition group (so $\#G_\infty = 1$ or $2$ according to whether $F$ is real or complex), then 
\[ U \cong (\Ind_{G_\infty}^G 1_{G_\infty}) - 1_G \cong \bigg\{ (x_\sigma) \in \bigoplus_{\sigma \in G/G_\infty} \!\!\! \Q \bigg\} / (1, 1, \dots, 1). \]
We call an element $u \in U$ corresponding to $(1,0, \dots, 0)$ under this isomorphism a ``Minkowski unit," where the first component corresponds to the trivial coset in $G/G_\infty$.  In other words, 
the elements $\sigma(u)$ for $\sigma \in G/G_\infty$ generate $U$ as $\Q$-vector space, and the only $\Q$-relation amongst these elements is that their sum vanishes.

Consider the matrix $M$ whose rows and columns are indexed by $\sigma, \tau \in G$, respectively, with entries given by $\log_p(\tau^{-1} \sigma(u))$,  where $u$ is a Minkowski unit.  It suffices to show that the $n \times n$ matrix $M$ has rank $r$, as its row space is equal to the image of $\reg_p$.  The matrix $M$ is equal to
\[ \sum_{\sigma \in G} \rho(\sigma) \log_p(\sigma(u)), \]
where $\rho$ denotes the regular representation of $G$ (represented in the standard basis given by the elements of $G$).  This fact is easily checked directly, and is familiar from Frobenius' study of the group determinant.

Decompose $\rho$ into irreducible representations $\rho_\pi$ over $\Q$, where $\pi$ denotes the representation space of $\rho_\pi$.  Note that the rational representations $\pi$ are not necessarily irreducible over $\overline{\Q}$, and they may appear with multiplicity.  We then have:
\[ \rank_{\C_p} M = \sum_\pi \rank_{\C_p} M_\pi, \qquad M_\pi = \sum_{\sigma  \in G} \rho_\pi(\sigma) \log_p \sigma(u). \]

Note that $M_\pi = 0$ if $\pi$ is the trivial representation.  Suppose now that $\pi$ is not trivial and let $\chi_\pi$ denote the character of $\pi$.  Let $\pi^+ = \pi^{G_\infty}$ denote the subspace on which $G_\infty$ acts trivially.

We will prove that
\[ \rank_{\C_p} M_\pi = \dim_{\Q} \pi^+ =  \frac{1}{\#G_\infty} \sum_{c \in G_\infty} \chi_\pi(c).\]
Define  $M_{\pi^+} \in \End_{\C_p}(\pi^+ \otimes_\Q \C_p)$ as the composition 
\begin{equation}\begin{tikzcd}
\pi^+ \otimes_\Q \C_p \ar[r,hook] &  \pi \otimes_\Q \C_p  \ar[r,"M_\pi"] & \pi \otimes_\Q \C_p \ar[r,"\pr^+"] &
  \pi^+ \otimes_\Q \C_p, \end{tikzcd}
  \end{equation} where \[ \pr^+ = \frac{1}{\#G_\infty}\sum_{c \in G_\infty} \rho_\pi(c) \] denotes the canonical projection 
  $\pi \rightarrow \pi^+$.

We claim that $M_{\pi^+}$ is invertible.  Suppose that this is not the case.
Then by Conjecture~\ref{c:mainq}, there exist nonzero rational vectors $w \in (\pi^+)^*, v \in \pi^+$ such that 
\[ \sum_{\sigma \in G} \langle w, \rho_\pi^+(\sigma)v \rangle \log_p \sigma(u) =0, \]
where $\rho_{\pi}^+ = \pr^+ \circ \rho_\pi$. Since $v \in \pi^+$ and $G_\infty$ fixes $u$ this yields
\[ \sum_{\sigma \in G/G_\infty} \langle w, \rho_\pi^+(\sigma)v \rangle \log_p \sigma(u) =0. \]

The only relation over $\Q$ satisfied by the $\log_p \sigma(u)$ for $\sigma \in G/G_\infty$ is that their sum vanishes.
 We therefore find that
$\langle w, \rho_\pi^+(\sigma)v \rangle$ is constant as $\sigma$ ranges over $G$.  Extend $w$ to an element $\tilde{w} \in \pi^*$ by defining $\tilde{w}(x) = w(\pr^+(x))$.  We  find that  $\langle \tilde{w}, \rho_\pi(\sigma)v \rangle$ is constant for $\sigma \in G$.
 This implies that $\tilde{w}$ is orthogonal to $(\rho_\pi(\sigma) - \rho_\pi(\tau))v$ for all $\sigma, \tau \in G$.  But \[ V = \Span\{ (\rho_\pi(\sigma) - \rho_\pi(\tau))v: \sigma, \tau \in G \} \] is clearly a $G$-stable rational subspace of  $\pi$, hence by irreducibility $V=0$ or  $V=\pi$.  If $V=\pi$, then $\tilde{w} =0$, a contradiction.  If $V=0$ and $v \neq 0$, then the span of $v$ is a copy of the trivial representation inside $\pi$, which is a contradiction since $\pi$ is non-trivial and irreducible.  We therefore conclude that $v = 0$, another contradiction.  This proves the claim that $M_{\pi^+}$ is invertible.

Therefore
\begin{equation} \label{e:mpi}
 \rank_{\C_p} M_\pi \ge \rank_{\C_p} M_{\pi^+} =  \frac{1}{\#G_\infty}\sum_{c \in G_\infty} \chi_\pi(c) = \dim_{\Q} \pi^+. 
 \end{equation}
Summing over $\pi$ and writing $\chi$ for the character of the regular representation $\rho$, we obtain
\[   \rank_{\C_p} M \ge \left(\frac{1}{\#G_\infty}\sum_{c \in G_\infty} \chi(c)\right) - 1 = \frac{n}{\#G_\infty} - 1 = r, \]
where the $-1$ comes from the exception of the trivial representation (which occurs with multiplicity exactly 1 in $\rho$).
Since of course $\rank_{\C_p} M  \le \dim_\Q U \le r$,  this inequality must be an equality (as must have been the inequality in (\ref{e:mpi})).  This concludes the proof.
\end{proof}

The exact same method of proof yields the following.

\begin{theorem} Conjecture~\ref{c:mainq} implies the Gross--Kuz'min conjecture. \label{t:gk}
\end{theorem}

\begin{proof}  Again we may consider a finite Galois CM extension $F/\Q.$  Let $F^+$ denote the maximal totally real subfield of $F$.  We may assume that every prime ideal of $F^+$ above $p$ splits completely in $F$, or else there is nothing to prove.
Write $G = \Gal(F/\Q)$.  
Let \[ U = \cO_F[1/p]^*  \otimes_{\Z} \Q, \qquad U^- = U/(c+1)U, \] where $c \in G$ is the unique complex conjugation.
Fix a place $\fp$ of $F$ above $p$ and let $G_\fp \subset G$ denote its decomposition group.
Write 
\[ n = \#G = [F:\Q], \qquad r = \dim_{\Q} U^- = \frac{n}{2 \# G_\fp} = \# \text{ of primes of } F^+ \text{ above } p.
 \]
The Gross--Kuz'min conjecture is that Gross's regulator map 
\[ \gr_p \colon U^- \otimes_\Q \C_p \longrightarrow \bigoplus_{\tau \in G} \C_p, \qquad u \mapsto (\log_p \N_{F_\fp/\Q_p}\tau^{-1}(u)) \]
is injective, i.e.\ has image of dimension $r$.  The map $\gr_p$ is $G$-equivariant, and we can describe the $G$-module structure of $U^-$ explicitly:
\[ U^- \cong (\Ind_{G_\fp}^G 1_{G_\fp})/(c+1) \cong \bigg\{ (x_\sigma) \in \bigoplus_{\sigma \in G/G_\fp} \!\!\! \Q \bigg\} / (c+1). \]
We consider a unit $u \in U$ corresponding to $(1,0, \dots, 0)$ under this isomorphism, where the first component correponds to the trivial coset in $G/G_\fp$.  In other words, 
the elements $\sigma(u)$ for $\sigma \in G/G_\fp$ generate $U$ as $\Q$-vector space, and the $\Q$-relations amongst these elements are generated by $c \sigma(u) + \sigma(u) = 0$. 

Consider the matrix $M$ whose rows and columns are indexed by $\sigma, \tau \in G$, respectively, with entries given by $\log_p(\N_{F_\fp/\Q_p}\tau^{-1} \sigma(u))$,  where $u$ is a unit as above.  It suffices to show that the $n \times n$ matrix $M$ has rank $r$, as its row space is equal to the image of $\gr_p$.  The matrix $M$ is equal to
\[ \sum_{\sigma \in G} \rho(\sigma) \log_p(\N_{F_\fp/\Q_p}\sigma(u)), \]
where $\rho$ denotes the regular representation of $G$. 

Decompose $\rho$ into irreducible representations $\rho_\pi$ over $\Q$, writing $\pi$ for the representation space of $\rho_\pi$.  We then have:
\[ \rank_{\C_p} M = \sum_\pi \rank_{\C_p} M_\pi, \qquad M_\pi = \sum_{\sigma  \in G} \rho_\pi(\sigma) \log_p \N_{F_\fp/\Q_p}\sigma(u). \]

Note that since $c$ is the unique complex conjugation in $G$, it is central.  Thus $c$ must act as a scalar $\pm 1$ on each irreducible representation $\pi$.  If $\rho_\pi(c)$ is the identity then since $c(u) = -u$, it is clear that $M_\pi = 0$.  Therefore we need only consider $\pi$ such that $c$ acts as $-1$.
Fix such a $\pi$ and let $\chi_\pi$ denote the character of $\pi$.  Let $\pi^+ = \pi^{G_\fp}$ denote the subspace on which $G_\fp$ acts trivially.
We will prove that
\[ \rank_{\C_p} M_\pi = \dim_{\Q} \pi^+ =  \frac{1}{\#G_\fp} \sum_{c \in G_\fp} \chi_\pi(c).\]
Define  $M_{\pi^+} \in \End_{\C_p}(\pi^+ \otimes_\Q \C_p)$ as the composition 
\begin{equation}\begin{tikzcd}
\pi^+ \otimes_\Q \C_p \ar[r,hook] &  \pi \otimes_\Q \C_p  \ar[r,"M_\pi"] & \pi \otimes_\Q \C_p \ar[r,"\pr^+"] &
  \pi^+ \otimes_\Q \C_p, \end{tikzcd}
  \end{equation} where \[ \pr^+ = \frac{1}{\#G_\fp}\sum_{c \in G_\fp} \rho_\pi(c) \] denotes the canonical projection 
  $\pi \rightarrow \pi^+$.

We claim that $M_{\pi^+}$ is invertible.  Suppose that this is not the case.
Then by Conjecture~\ref{c:mainq}, there exist nonzero rational vectors $w \in (\pi^+)^*, v \in \pi^+$ such that 
\[ \sum_{\sigma \in G} \langle w, \rho_\pi^+(\sigma)v \rangle \log_p \N_{F_\fp/\Q_p}\sigma(u) =0, \]
where $\rho_{\pi}^+ = \pr^+ \circ \rho_\pi$. Since $v \in \pi^+$ and $G_\fp$ fixes $u$, while $c$ acts as $-1$ on $v$ and $u$, this yields
\[ \sum_{\sigma \in S} \langle w, \rho_\pi^+(\sigma)v \rangle \log_p \N_{F_\fp/\Q_p} \sigma(u) =0 \]
where $S$ is a fixed set of representatives for $G/\langle G_\fp, c \rangle$.
Extend $w$ to an element  $\tilde{w} \in \pi^*$ by defining $\tilde{w}(x) = w(\pr^+(x))$, so we can write the previous equation as
\[ \sum_{\sigma \in S} \langle \tilde{w}, \rho_\pi(\sigma)v \rangle \log_p \N_{F_\fp/\Q_p} \sigma(u) =0. \]
Scale so that $u \in \cO_F[1/p]^*$ and such that the scalars $\langle \tilde{w}, \rho_\pi(\sigma)v \rangle$  are integers.  Since
\[ \prod_{\sigma \in S} \N_{F_\fp/\Q_p} \sigma(u)^{ \langle w, \rho_\pi(\sigma)v \rangle} = 
\prod_{\sigma \in S} \prod_{g \in G_\fp} g\sigma(u)^{ \langle w, \rho_\pi(\sigma)v \rangle} \]
lies in the kernel of $\log_p$, it is a power of $p$ times a root of unity.  Since all the conjugates of $u$ have complex absolute value 1, the power of $p$ must be trivial and the product is a root of unity.  (Of course upon tensoring with $\Q$, roots of unity vanish.)  Employing the change of variables $g \sigma \sim \tau$, where $\sim$ denotes equivalence mod 
$\langle G_\fp, c \rangle$, we obtain that
\[  \prod_{\tau \in S} \tau(u)^{ \langle w, \sum_{g \in G_\fp} \rho_\pi(g^{-1}\tau)v \rangle}
\]
is trivial.
Yet by definition of $u$, the $\tau(u)$ for $\tau \in S$ are linearly independent over $\Q$, so we obtain 
\[ \langle w, \sum_{g \in G_\fp} \rho_\pi(g^{-1}\tau)v \rangle = 0 \]
for all $\tau$.  By construction $w \circ \rho_\pi(g) = w$ for all $g \in G_\fp$, so we simply have
\begin{equation} \label{e:worth}
  \langle w,  \rho_\pi(\tau)v \rangle = 0 \end{equation}
   for all $\tau \in S$.
Now $v$ is stable under $G_\fp$ and is negated by $c$, so this equation holds for all $\tau \in G$.  By the irreducibility of $\pi$ and the fact that $v \neq 0$, the span of the $\rho_\pi(\tau)v$ for all $\tau \in G$ is all of $\pi$.  Equation (\ref{e:worth}) then contradicts $w \neq 0$.
 This proves the claim that $M_{\pi^+}$ is invertible.

We then proceed as in the proof of Theorem~\ref{t:leopoldt}. We have
\begin{equation} \label{e:rankmpi1}
 \rank_{\C_p} M_\pi \ge \rank_{\C_p} M_{\pi^+} =  \frac{1}{\#G_\fp}\sum_{g \in G_\fp} \chi_\pi(g) = \dim_{\Q} \pi^+. 
 \end{equation}
 To write an inequality and equation that holds for all irreducible $\pi$ (i.e.\ including those on which $c$ acts as $+1$) we write
 \begin{equation}  \label{e:rankmpi2}
 \rank_{\C_p} M_\pi \ge \rank_{\C_p} M_{\pi^+} =  \frac{1}{\#G_\fp}\sum_{g \in G_\fp} \frac{\chi_\pi(g) - \chi_\pi(gc)}{2}.
 \end{equation}
 
Summing over $\pi$ and writing $\chi$ for the character of the regular representation $\rho$, we obtain
\[   \rank_{\C_p} M \ge \left(\frac{1}{\#G_\fp}\sum_{g \in G_\fp} \frac{\chi(g) - \chi(gc)}{2} \right)  = \frac{n}{2\#G_\fp} = r. \]
Note that in the first equality we use that $c \not \in G_\fp$.
Since of course $\rank_{\C_p} M  \le \dim_\Q U \le r$,  this inequality must be an equality (and the same is therefore true for the inequalities (\ref{e:rankmpi1}) and (\ref{e:rankmpi2})).  This concludes the proof.
\end{proof}

\section{The Theorems of Masser and Waldschmidt} \label{s:mw}

To motivate our approach to the Matrix Coefficient Conjecture, we recall simplified versions of the theorems of Masser \cite{masser} and Waldschmidt \cite{w}, as relevant for our present context.  We first set notations. Let $F$ be a number field, and let \[ x_i = (x_{ij})_{j=1}^{n} \in (F^*)^n, \qquad i = 1, \dotsc, m. \]   Let $X \subset (F^*)^n$ be the subgroup generated by the $x_i$.  There is a bilinear pairing \[ \langle \ \! , \rangle_X : \Z^m \times \Z^n \rightarrow F^* \] given by
\[ \langle (a_1, \dotsc, a_m), (b_1, \dotsc, b_n) \rangle_X = \prod_{i,j} x_{i,j}^{a_i b_j}. \]  
Consider the following condition, which we call ``Condition O", for orthogonal:
\begin{quote}
(O) There exist subgroups $A \subset \Z^m$ and $B \subset \Z^n$ of ranks $m', n'$, respectively, with $\langle A, B \rangle_X = 1$ and $m'/m + n'/n > 1.$
\end{quote}

For a positive real $N$, we let
\begin{align*}
 \Z^m(N) &= \{ a \in \Z^n : 0 \le a_i \le N-1, i=1, \dotsc, m \} \subset \Z^m, \\
 X(N) &= \{ x_1^{a_1}x_2^{a_2} \cdots x_m^{a_m} : a \in \Z^m(N) \} \subset X. 
\end{align*}
Consider the following condition, which we call ``Condition M", for Masser:
\begin{quote}
(M) There exists a polynomial $P(t_1, \dotsc, t_n) \in F[t_1, \dotsc t_n]$ of total degree $< N^{m/n}$ such that $P(x) = 0$ for all 
$x \in X(nN)$.
\end{quote}

Masser's Theorem \cite{masser}*{Theorem 2} is the following:
\begin{theorem}[Masser]
If Condition $M$ holds for any positive integer $N$, then Condition O holds.
\end{theorem}

On the transcendence side, we have the result of Waldschmidt. 
 Consider the square matrix \[ A = (\log(x_{ij})) \in \C^{n \times n} \qquad \text{or} \qquad 
 A = (\log_p(x_{ij})) \in \C_p^{n \times n}.
 \]
(In the second case,  we suppose an embedding $F \subset \C_p$ has been fixed, and we assume that the $x_{ij}$ are all $p$-adic units.)
Consider the following condition, which we call ``Condition W", for Waldschmidt:
\begin{quote}
(W) The rank of $A$ is less than $mn/(m+n)$. 
\end{quote}

\begin{theorem}[Waldschmidt]  If Condition W holds, then Condition $M$ holds for all sufficiently large $N$.
\end{theorem}

Combining these results, one obtains the celebrated theorem of Waldschmidt--Masser, which states $\text{(W)} \Longrightarrow \text{(O)}$.
As a corollary, one obtains that the ranks of the Leopoldt matrix and of the Gross regulator matrix are at least half their expected values,
 since in these contexts, condition (O) leads to a contradiction (see \cite{das}*{Corollary 5.3}).  
 
  It is also worth noting that (in the general case) the converse implication to Masser's Theorem,
(O) $\Longrightarrow$ (M), is  not hard to prove. 

\section{Proof Strategy} \label{s:strategy}

With the goal of attacking Conjecture~\ref{c:mainq}, we propose the following modifications of the conditions (O), (W), and (M) from \S\ref{s:mw}. We are now working with a square matrix ($m=n$), but the rest of the notation is unchanged.

\begin{quote}
(o) There are nonzero elements $a, b \in \Z^n$ such that $\langle a, b \rangle_X = 1$.
\end{quote}
\begin{quote}
(w) The matrix $A$ has rank less than $n$.
\end{quote}
For the last condition we introduce some notation. Throughout this section we let $F$ denote $\C$ or $\C_p$ (or more generally, any algebraically closed field of characteristic 0).
For a polynomial $P \in F[t_1, \dots, t_n]$, we denote by $S(P)$ the support of $f$, i.e.\ the set of tuples $(a_1, \dots, a_n) \in \Z^n$ such that the coefficient of $t_1^{a_1} \cdots t_n^{a_n}$ in $P$ is nonzero.  Let $\Delta(P) \subset \R^n$ denote the convex hull of $S(P)$.

\begin{definition} The Bernstein--Kushnirenko number of $P$ is defined to be \[ \BK(P) = n! \Vol(\Delta(P)). \]
More generally, given $f_1, \dots, f_n \in F[t_1, \dots, t_n]$, we define
\[ \BK(f_1, \dots, f_n) = n! \Vol(\Delta(f_1), \dots, \Delta(f_n)) \]
where the right side denotes mixed volume.  In particular, $\BK(P, P, \dots, P) = \BK(P)$.  
\end{definition}

The key to our condition (m) is the following remarkable theorem of Bernstein and Kushnirenko.
For details on the notion of mixed volume and the Bernstein--Kushnirenko theorem, see \cite{clo}*{Sections 7.4 and 7.5}.

\begin{theorem}[Bernstein--Kushnirenko] \label{t:bk}  Let $f_1, f_2, \dots, f_n \in F[t_1, \dots, t_n]$  have finitely many common zeroes in $(F^*)^n$.
Then  the number of common zeroes of the $f_i$ in $(F^*)^n$ is bounded by $\BK(f_1, \dots, f_n)$. 
\end{theorem}  
In fact the number of common zeroes in $(F^*)^n$ is precisely $\BK(f_1, f_2, \dots, f_n)$ under an appropriate genericity assumption, but we will not use this result.

\begin{definition} Let $P \in F[t_1, \dots, t_n]$.  Define the {\em BK-degree} of $P$, denoted $\BKd(P)$, to be the maximum as $i = 1, \dots, n$ of
$\BK(P, \dots, P, \ell, \dots, \ell)$, where $\ell$ denotes a generic linear polynomial and there are $i$ copes of $P$ and $(n-i)$ copies of $\ell$.
\end{definition}

We can now state our last condition on the matrix $M$.

\begin{quote}
(m) For  $N$ large enough, there exists a polynomial $P(t_1, \dotsc, t_n)$ such that  $\BKd(P) < N^n$ and $P(x)=0$  for all $x \in X(nN)$.
\end{quote}

\begin{remark}  If $P$ is not divisible by any monomial, then $\BK(P, \ell, \dots, \ell)$ is the usual degree of $P$.  We may assume that $P$ is not divisible by any monomial by dividing by such a monomial if it exists, since monomials have no zeroes in $(F^*)^n$.  Therefore, when $n = 2$, the condition $\BKd(P)  < N^2$ is equivalent to $\deg(P) < N^2$ and $\BK(P) < N^2$.
\end{remark}

\begin{remark}
The  condition $\BKd(P) < N^n$ is weaker than the counterpart $\deg(P) < N$ in the setting of Masser's Theorem. For example, $\deg(P) < N$ implies that $S(P)$ is contained  in the standard simplex of side length $N$, which has volume $(1/n!)N^n$.  Even when $n=2$, our conditions on the polynomial $P$ given in the remark above are strictly weaker than the condition $\deg(P) < N$.
\end{remark}

Our strategy for proving Conjecture~\ref{c:mainq} is to attempt to prove that
\[ \text{(w)} \Longrightarrow \text{(m)} \Longrightarrow \text{(o)}. \]
The following is not logically necessary, but provides reassurance that this is a reasonable strategy:

\begin{prop}  (o) $\Longrightarrow$ (m) \label{p:oim}
\end{prop}

\begin{proof}   Suppose that $w_0Av_0 = 0$, with $w_0, v_0 \in \Z^n$ nonzero. 
Write  $w_0 = (w_1, \dotsc, w_n)$  and consider the monomial \[ t_1^{w_1}t_2^{w_2} \cdots t_n^{w_n} \in \Z[t_{i}^{\pm 1}]. \]
 Evaluated at the $n$-tuple $(e^{e_1Av}, e^{e_2 Av}, \dots, e^{e_nAv})$, this equals $e^{w_0 A v}$.  Since $w_0 A v_0 = 0$, this expression only depends on the image of $v$ in $\Z^n/\langle v_0 \rangle$.  Therefore,  as $v$ ranges over the $(nN)^n$ elements of $\Z^n(nN)$,
 the number of distinct values $w_0 A v_0$ is bounded by $C \cdot N^{n-1}$ for some constant $C$ not depending on $N$. 
 
 Therefore, we define
\[ P_0(t) := \prod_{ \text{distinct values of } w_0Av \text{ as } v \in \Z^n(N)} (t_1^{w_1}t_2^{w_2} \cdots t_n^{w_n} - e^{w_0 A v}). \]

We multiply by a large enough monomial to clear away the negative exponents (but such that the result is not divisible by any monomial).  This yields a polynomial $P$ 
 such that $P(x) = 0$ for all $x \in X(nN)$.    The degree of $P$ is bounded by a constant times the number of terms in the product, which, as described above, is $O(N^{n-1})$.  Therefore, for $N$ sufficiently large, the degree of $P$ is less than $N^n$.  
 
 Now, the support $S(P)$ lies on a line segment.  Since $n > 1$, an elementary argument shows that
 \[ \BK(P, \dots, P, \ell, \dots, \ell) = \begin{cases}
 \deg(P) & i = 1 \\
0 & i > 1,
  \end{cases}
 \]
 where $i$ is the number of copies of $P$ in the argument.  The result follows.
\end{proof}

In view of Proposition~\ref{p:oim}, since (w) is expected to imply (o) by Theorem~\ref{t:imply}, we should indeed expect 
(w) to imply (m).  Unfortunately we do not know how to prove this statement.

\section{Proof of Analogue of Masser's Theorem}

In this section we prove:

\begin{theorem}  (m) $\Longrightarrow$ (o). 
\label{t:mo}
\end{theorem}

The idea of the proof follows Masser's, but uses the condition $\BKd(P) < N^n$ rather than a bound $\deg(P) < N$ to control certain intersections.
We let the group $X \subset (F^*)^n$ act on the Laurent polynomial ring $R = F[t_1^{\pm 1}, \dots, t_n^{\pm 1}]$ by \[ z\cdot f = f(z_1 t_1, z_2 t_2, \dots, z_nt_n). \]

 Suppose  we have a polynomial $P(t)$ as given by condition (m).   Given elements $y_i \in X((n-1)N)$ and constants $c_i \in F$, any linear combination $\sum c_i (y_i \cdot P)$ vanishes on $X(N)$.  Under the assumption that condition (o) does not hold, we will use the condition $\BKd(P) < N^n$ to show that there exist $n$ such linear combinations, say $f_1, \dots, f_n$, which have finitely many common zeroes.  Note that $S(f_i) \subset S(P)$ for each $i$.  Hence 
\[ \BK(f_1, f_2, \dots, f_n) \le \BK(P) < N^n. \]
Yet $X(N)$ has $N^n$ elements (unless (o) holds for some nonzero $a$ and {\em every} $b$), so we obtain a contradiction.  This implies that (o) must hold.

Let us now describe the existence of the $f_i$.
Recall that the subgroup $X \subset (F^*)^n$ is generated by elements $x_1, \dots, x_n$.
If $a \in \Z^n$, we write $x^a = \prod_{i=1}^{n} x_i^{a_i} \in X$. 
For a prime ideal $\fp \subset R$, let 
\[ \Stab_X(\fp) = \left\{a \in \Z^n \colon x^a \cdot \fp = \fp\right\}. \]
Consider the maximal ideal $\fm = (t_1 -1, \dots, t_n - 1)$.  The following lemma follows from \cite{das}*{Lemma 5.9}.
\begin{lemma} Suppose there exists a nonzero prime ideal $\fp \subset \fm$ such that $\Stab_X(\fp) \neq 0$. Then condition (o) holds.
\end{lemma} 
Therefore, we may assume that $\Stab_X(\fp) = 0$ for each prime ideal $\fp \subset \fm$.  Now, the maximal ideals of $R$ all have the form $z \cdot \fm$ for some $z \in (F^*)^n$.  Hence any prime ideal $\fp \subset R$ satisfies $\fp \subset z \cdot \fm$, i.e., $z^{-1} \cdot \fp \subset \fm$ for some $z$.  Since $\Stab_X(\fp) = \Stab_X(z^{-1} \cdot \fp)$, we may  assume that $\Stab_X(\fp) = 0 $ for all nonzero prime ideals $\fp$.

We may now construct the polynomials $f_i$ described above.  For $i = 1, \dots, n$, we inductively construct a nontrivial linear combination
\begin{equation} \label{e:fform}
 f_i = \sum_{j} c_j y_j \cdot P \end{equation}
with $c_j  \in F$, $y_j \in X((i-1)N)$ such that the zero set $V(f_1, \dots, f_i)$  has dimension $n - i$ and fewer than $N^n$ irreducible components in $(F^*)^n$.
For the base case $i=1$, we take $f_1 = P$.  The zero set $V(P)$ has codimension 1 and the number of irreducible components is equal to the number of distinct irreducible factors of $P$.   This is bounded by $\deg(P) = \BK(P, \ell, \ell, \dotsc, \ell)$, which by assumption is less than $N^n$.

For the inductive step, suppose we have constructed $f_1, \dots, f_i$ such that $V(f_1, \dots, f_i)$  has dimension $n - i$ and fewer than $N^n$ irreducible components.  Each irreducible component corresponds to a nonzero prime ideal $\fp \subset R$.  Each such prime satisfies $\Stab_X(\fp) = 0$ and therefore the $N^n$ primes $z^{-1} \cdot \fp$ for $z \in X(N)$ are distinct.  (As mentioned above, $|X(N)|= N^n$ unless (o) holds.)  Since $V(f_1, \dots, f_i)$ has fewer than $N^n$ irreducible components, there exists at least one $z \in X(N)$ such that $V(z^{-1} \cdot \fp)$ is not an irreducible component of $V(f_1, \dots, f_i)$.  Therefore, for each such $\fp$, there exists an $f_j$ for $j = 1, \dots, i$ such that $f_j \not \in z^{-1} \cdot \fp$, i.e., $z \cdot f_j \not \in \fp$. 

By \cite{das}*{Lemma 5.13}, there exists an $F$-linear combination of the $z \cdot f_j$ that does not lie in {\em any} such $\fp$.  We define $f_{i+1}$ to be this linear combination.
We check that it has the desired properties.  Firstly, since each $f_1, \dots, f_i$ has the form (\ref{e:fform}) with $y_j \in  X((i-1)N)$, it follows that $f_{i+i}$ has the same form with $y_j \in X(iN)$.  Moreover, by construction, $f_{i+1} \not \in \fp$ for the prime $\fp$ corresponding to any irreducible component of $V(f_1, \dots, f_i)$, so it follows that $V(f_1, \dots, f_{i+1})$ has dimension $n - (i+1)$.  Finally, the number of irreducible components of $V(f_1, \dots, f_{i+1})$ is bounded by its degree, which is in turn bounded by  
$\BK(f_1, \dots, f_{i+1},\ell, \dots, \ell)$ by Theorem~\ref{t:bk}.  Since each $\Delta(f_i) \subset \Delta(P)$, this is bounded by $\BK(P, \dots, P,\ell, \dots, \ell)$, which by assumption is less than $N^n$.  This completes the inductive construction of the $f_i$.

To complete the proof, we note that in the final stage, we have constructed $f_{1}, \dots, f_n$ such that $V(f_{1}, \dots, f_n)$ is a collection of fewer than $N^n$ points.  However, as noted at the outset, each $f_i$ vanishes on $X(N)$, which has size $N^n$.  This is a contradiction.  Since the only assumption made was that condition (o) does not hold, the proof of Theorem~\ref{t:mo} is complete.

\section{Discussion}

\subsection{Fewnomials}

A proxy for $\BK(P)$ is the size of the support, $|S(P)|$.  Indeed, if $S(P)$ contains all the lattice points in 
$\Delta(P)$, then $\BK(P) \approx n! |S(P)|$.  For $n=2$, an exact relationship is provided by Pick's theorem, which distinguishes between points on the boundary and interior of $\Delta(P)$.
Let us consider the following condition, which is less stringent than condition (m):
\begin{quote}
(m') For  $N$ large enough, there exists a polynomial $P(t_1, \dotsc, t_n)$ such that  $|S(P)| < N^n$ and $P(x)=0$  for all $x \in X(2N)$.
\end{quote}

The main result of this section is the following.

\begin{theorem} \label{t:wmo}
  Suppose that $A \in M_2(\R)$ or $M_2(\C_p)$ satisfies condition (w).
 In the $p$-adic case, suppose we can write $A = (a_i b_j)$ for two vectors $a, b \in \C_p^2$ with $|a_i|_p \le p^{-p/(p-1)}$ and $|b_i|_p \le 1$, which can be obtained by scaling $A$ by a power of $p$.
   Then  (m') $\Rightarrow$ (o) for the matrix $A$.
\end{theorem}

The crux of the proof is a bound on the number of roots of a function that is a sum of exponentials.
In the real case, the result is elementary and follows from \cite{lil}*{Lemma 2.1}.
\begin{lemma} \label{l:rootsr} (real case).  Let $n \ge 1$ and let $f(z) =\sum_{i=1}^{n} b_i e^{w_i z}$ with $b_i, w_i \in \R$.
Then $f(z)$ has at most $n-1$ zeroes.
\end{lemma}

The $p$-adic case is slightly more involved.

\begin{lemma} \label{l:roots} ($p$-adic case).   Let $n \ge 1$ and let $f(z) =\sum_{i=1}^{n} b_i e^{w_i z}$ with
$b_i, w_i \in \C_p$ such that  $|w_i | \le p^{-p/(p-1)}$. Then $f(z)$ has at most $\frac{p}{p-1}(n-1)$ zeroes in the disc $|z | \le 1$.
\end{lemma}

\begin{proof}[Proof of Lemma~\ref{l:roots}] 
We use the $p$-adic Schwarz lemma (see \cite{das}*{Lemma 5.7}), which implies that for any $R > 1$, the number of roots of $f(z)$ in the disc $|z| \le 1$ is bounded as follows:
\begin{equation} 
\label{e:nr}
 \text{ number of roots of } f(z) \text{ in the unit disc } \le \frac{\log (|f|_R/|f|_1)}{\log R}. 
 \end{equation} 

We will fix $R = p.$  We want to bound $|f|_R/|f|_1$.
For any  $u \in \C_p$, there is a unique polynomial $P_u(z) = \sum_{i=1}^{n} a_i(u) z^{i-1}$ of degree at most $n-1$ such that $P_u(w_i) = e^{w_i u}$ for $i = 1, \dotsc, n$.  An elementary computation then shows that \begin{equation} f(u) = \sum_{i=1}^{n} a_i(u) f^{(i-1)}(0). \label{e:fsum} \end{equation}

Now \[ |f|_1 = \sup_{i \ge 0} | \frac{f^{(i)}(0)}{i!}| \ge \sup_{i \ge 0} | f^{(i)}(0)|, \] so it follows from (\ref{e:fsum}) that
\begin{equation} |f|_R \le \sup_{i=1}^{n} |a_i|_R \cdot |f|_1.  \label{e:frat} \end{equation}
For simplicity write $x = p/(p-1)$.  It remains to prove that $|a_i|_R \le p^{x(n-1)}$.   Combining (\ref{e:nr}) and (\ref{e:frat}) then gives the desired result.
Now one can prove by induction that
\[ P_u(z) = \sum_{i=1}^{n} c_i(u) \prod_{j = 1}^{i-1} (z - w_j) \]
where \[ c_i(u) = \sum_{s=1}^{i} \frac{e^{w_s u}}{\prod_{s' \le i, s' \neq s} (w_s - w_{s'})}. \]
(Note that this is just writing out (6.3.5) and (6.3.6) in \cite{lil} explicitly in our case.)  
Clearly it suffices to show that $|c_i|_R  \le p^{x(n-1)}$, as the same bound will then hold for $|a_i|_R$.  In fact we will show that $|c_i|_R  \le p^{x(i - 1)} \le p^{x(n-1)}$.
Viewing $i$ as fixed for the moment, we define for  any nonnegative integer $k$: 
\[ d_k = \sum_{s=1}^{i} \frac{w_s^k}{\prod_{s' \le i, s' \neq s} (w_s - w_{s'})}, \qquad \text{ so } c_i(u) = \sum_{k=0}^{\infty} \frac{d_k u^k}{k!}. \]
We will show that $|d_k| \le p^{-x(k - i + 1)}.$  Indeed, $d_k$ is given explicitly as the coefficient of $z^{i-1}$ when the polynomial $z^k$ is reduced modulo $\prod_{j=1}^{i} (z - w_i)$.  In particular,
\[ d_0 = d_1= \dotsc = d_{i-2}  =0, \  d_{i-1} = 1, \] and
\[ d_k = s_1 d_{k-1} + s_2 d_{k-2} + \cdots + s_i d_{k-i}, \qquad k \ge i \]
where $s_j$ are the symmetric polynomials in $w_1, \dotsc, w_i$ with appropriate sign.  The desired bound on the $d_k$ then follows by induction since $|w_i| \le p^{-x}$.  We thus find that
\[ |c_i|_R \le \sup_{k\ge 0} \frac{p^{-x(k - i + 1)} p^{k}}{| k! | } = p^{x(i - 1)} \cdot  \sup_{k\ge 0} \frac{p^{(1-x)k}}{| k! | } \le p^{x(i - 1)}. 
\]
This concludes the proof.
\end{proof}

We can now prove Theorem~\ref{t:wmo}.

\begin{proof}[Proof of Theorem~\ref{t:wmo}]
We  write $A = (a_i b_j)$ for two vectors $a, b \in F^2$, where $F = \R$ or $\C_p$. In the $p$-adic case we assume the conditions on $a_i, b_i$ as in the statement of the Theorem.
Then $wAv$ simply becomes $(w \cdot a)(v \cdot b)$.   The polynomial $P(t) = \sum_{w} c_w t^w$ yields a function $\sum_w c_w e^{(w \cdot a)z}$ that has the roots $(v \cdot b)$ for $v \in \Z(2N)$.  Since there are $4N^2$ such $v$'s and there are  fewer than $N^2$ $w$'s in the sum, this is a contradiction to Lemma~\ref{l:rootsr} or~\ref{l:roots}, unless two of the $w \cdot a$'s are equal or two of the $v \cdot b$'s are equal.  This gives the desired result that the $a$'s or $b$'s are linearly dependent over $\Z$, which is exactly what condition (o) says in this case.
\end{proof}

Theorem~\ref{t:wmo} implies that, in order to prove the 4 exponentials conjecture over $\R$ or $\C_p$, it suffices to prove that (w) $ \Rightarrow$ (m') instead of the stronger implication (w) $\Rightarrow$ (m).  Of course, this implication remains a mystery.

\subsection{Auxiliary Polynomials}

We conclude by pointing out the difficulty in constructing the auxliary polynomials necessary in (m).  As an instructive example, consider the matrix
\[ A = \begin{pmatrix}
 0 & -\log \alpha & -\log \beta  \\
\log \alpha & 0 & -\log \gamma &  \\
 \log \beta & \log \gamma & 0 & 
\end{pmatrix}
 \]
which is skew-symmetric and hence has vanishing determinant for any $\alpha, \beta, \gamma$.  The polynomials in condition (m) that we are trying to construct (witnessing the 0's in the matrix $A$) are 
\begin{equation} \label{e:pdef}
P(t_1, t_2, t_3)= \prod_{i,j = 0}^{3N-1} (t_1 - \alpha^i \beta^j) 
\end{equation}
and similar poynomials in $t_2, t_3$.  For simplicity suppose that $\alpha, \beta, \gamma$ are integers $\ge 2$, and let $M$ be their maximum.  Then the coefficients of $P$ have size on the order of $M^{O(N^3)}$, so a naive bound on $|P(x,y,z)|_\infty$ for any arguments $x,y, z$ will be at least this large.  Let us assume we are working with the $p$-adic conjecture, and that the integers $\alpha, \beta, \gamma$ are congruent to 1 modulo $p$.  Then the $9N^2$ roots of $P$ (as a function of $t_1$) imply that $|P(x, y, z)|_p < p^{-9N^2}$ for any $x \in 1 + p \Z_p$ (in particular for any $(x,y,z) \in X$).   The usual method in ($p$-adic) transcendence to prove that a polynomial vanishes at an algebraic (say, integer) argument $(x,y,z)$ is to show that 
\[ |P(x,y,z)|_\infty \cdot |P(x,y,z)|_p < 1. \]
However, the naive bounds we have given above for our polynomial $P$ do not allow for such a deduction.  One might say that the reason that the polynomial $P$ vanishes on $X(3N)$ is not ``analytic.'' 
This leads to the following question, to which we do not know the answer (even for the particular matrix $A$ given above):
\begin{question}  For every $N$ large enough, does there exist a polynomial $P \in \Z[t_1, t_2, t_3]$ and functions $f(N), g(N)$ satisfying the following?
\begin{itemize}
\item $\BKd(P) < N^3$.
\item $|P(x,y,z)|_p < f(N)$ for all $(x,y,z) \in X$.
\item $|P(x,y,z)|_\infty < g(N)$ for all $(x,y,z)$ with $|x|, |y|, |z| < M^{3N}.$
\item $f(N) \cdot g(N) < 1$ for $N$ sufficiently large.
\end{itemize}
\end{question}

We say that such a polynomial $P$ vanishes on $X(3N)$ for analytic reasons.  The construction of such polynomials (such as in the proof of Waldschmidt's theorem) fits into the usual framework of transcendence theory.  Constructing the polynomial in (\ref{e:pdef}) seems to require a different approach.

\end{document}